\documentclass[12pt]{article}
\usepackage{amsmath}
\usepackage{amssymb}
\usepackage{graphicx}
\usepackage[cp1251]{inputenc}
\usepackage[russian]{babel}
\newcommand{\il}[2]{\int\limits_{#1}^{#2}}
\newcommand{\ilp}[1]{\int\limits_{#1}^{+\infty}}

\newcommand{\ph}{\phantom{a}}
\newcommand{\phh}{\phantom{aaa}}

\newcommand{\sist}[2]{\left\{
\begin{array}{l}
{#1}\\
\ph\\
{#2}
\end{array}
\right.}

\begin{document}

\vskip 20pt

MSC 34C10

\vskip 10pt

\centerline{\bf Oscillation, suboscillation and nonoscillation criteria for }
 \centerline{\bf linear  systems of ordinary differential equations}

\vskip 5 pt

\centerline{\bf G. A. Grigorian}
\vskip 10 pt

\noindent
Abstract. The Riccati equation method and an approach of the use of unknown factors is used to establish oscillation, suboscillation and nonoscillation criteria for linear systems of ordinary differential equations. A necessary condition for Lyapunov (asymptotic) stability for these systems is obtained.

\vskip 5 pt


\vskip 5 pt

{\bf 1. Introduction}. Let $a_{jk}(t), \ph j,k = \overline{1,n}$ be real-valued locally integrable functions on $[t_0,\infty)$. Consider the linear system of ordinary differential equations
$$
\phi_j' = \sum\limits_{k=1}^n a_{jk}(t) \phi_k, \ph j=\overline{1,n}, \ph n\ge 3, \ph t \ge t_0.  \eqno (1.1)
$$

{\bf Definition 1.1} {\it The system (1.1) is called oscillatory, if for its every solution \linebreak $(\phi_1(t),\dots,\phi_n(t))$ the function $\phi_1(t)$ has arbitrarily large zeroes.
}

{\bf Definition 1.2} {\it The system (1.1) is called suboscillatory, if for its every solution \linebreak $(\phi_1(t),\dots,\phi_n(t))$  at least one of the functions $\phi_1(t), \ph \phi_3(t),\dots,\phi_n(t)$ has arbitrarily large zeroes.
}

{\bf Definition 1.3} {\it The system (1.1) is called nonoscillatory, if  it has a  solution \linebreak $(\phi_1(t),\dots,\phi_n(t))$  such that  $\phi_k(t) \ne 0, \ph t \ge T$ for some $T\ge t_0, \ph k=\overline{1,n}$.
}

{\bf Remark 1.1.} {\it If the system (1.1) is nonoscillatory, the it is neither oscillatory nor suboscillatory.}

The questions of oscillation and nonoscillation of linear systems of ordinary differential equations are important problems of qualitative theory of differential equations. These problems are studied by many authors specially in the case of Hamiltonian systems due to their important applications in the natural sciences (see [1,2,6,7,11-16] and cited works therein). However it should be noticed here that the definitions of oscillation and nonoscillation for Hamiltonian systems are radically different from the Definitions 1.1 and 1.3. Some oscillation and nonoscillation results  in the particular cases $n =2$ and  $n =4$ for the system (1.1) have been obtained in the papers [4,5,8.10] (see also [7]).

In this paper the Riccati equation method a new  approach (which we call the unknown factors approach) is used to establish oscillation, suboscillation and nonoscillation criteria for the system (1.1).  A necessary condition for Lyapunov (asymptotic) stability for this system is obtained.

{\bf 2. Auxiliary propositions}. Let $f_k(t), \ph g_k(t)$ and $h_k(t), \ph k=1,2$ be real-valued  locally integrable functions on $[t_0,+\infty)$. Consider the scalar Riccati equations
$$
y' + f_k(t) y^2 + g_k(t) y + h_k(t) = 0, \ph t\ge t_0, \ph k=1,2 \eqno (2.1_k)
$$
and the differential inequalities
$$
\eta' + f_k(t) \eta^2 + g_k(t) \eta + h_k(t) = 0, \ph t\ge t_0, \ph  k=1,2. \eqno (2.2_k)
$$

{\bf Remark 2.1.} {\it Every solution of Eq. $(3.1_k)$ on $[t_1,t_2) \ph (t_0 \le t_1< t_2 \le +\infty)$ is also a solution of the inequality $(3.2_k), \ph k =1,2.$}

{\bf Remark 2.2.} {\it If $f_k(t) \ge 0, \ph t\ge t_0,$ then every solution of the linear equation
$$
\zeta' + g_k(t) \zeta + h_k(t) = 0, \phh t\ge t_0
$$
is also a solution of the inequality $(3.2_k), \ph k =1,2.$}

The following comparison theorem is important in thee proof of the main results.

\vskip 5pt

{\bf Theorem 2.1.} {\it Let Eq. $(2.1_2)$ have a real-valued solution $y_2(t)$ on $[t_0,\tau_0) \ph (t_0 < \tau_0 \le +\infty)$ and let the following conditions be satisfied: $f_1(t) \ge 0$  and $\il{t_0}{t} \exp\biggl\{\il{t_0}{\tau}[f(s)(\eta_1(s) + \eta_2(s)) + g(s)]d s\biggr\}[(f_1(\tau) - f(\tau))y_2^2(\tau) + (g_1(\tau) - g(\tau)) y_2(\tau) + h_1(\tau) - h(\tau)]d \tau \ge 0, \ph t\in [t_0,\tau_0)$ where $\eta_1(t)$ and $\eta_2(t)$ are solutions of the inequalities $(2.2_1)$ and $(2.2_2)$ respectively on $[t_0, \tau_0)$ such that $\eta_j(t_0) \ge 2 y_2(t_0), \ph j=1,2.$ Then for every $\gamma_0 \ge y_2(t_0)$ Eq. $(2.1_1)$ has a solution $y_1(t)$ on $[t_0, \tau_0)$, satisfying the condition $y_1(t_0) = \gamma_0.$}

Proof. By analogy with the proof of Theorem 3.1 from [9].

\vskip 5pt

{\bf Remark 2.3.} {\it One can easily verify, that in the case $\tau_0 <+\infty$ Theorem 2.1 remains valid if we replace $[t_0,\tau_0)$ by $[t_0,\tau_0]$ in it.}

Let $p_{jk}(t), \ph j,k =1,2$ be real-valued locally integrable functions on $[t_0,+\infty)$. Consider the linear system of ordinary differential equations
$$
\sist{\phi' = p_{11}(t) \phi + p_{12}(t) \psi,}{\psi' = p_{21}(t) \phi + p_{22}(t) \psi, \ph t\ge t_0.} \eqno (2.3)
$$

{\bf Definition 2.1.} {\it The system (2.3) is called oscillatory if for its every solution  $(\phi(t), \psi(t))$ the function $\phi(t)$ has arbitrarily  large zeroes.}

{\bf Definition 2.2.} {\it The system (2.3) is called oscillatory  on the interval $[a,b]$  if for its every solution  $(\phi(t), \psi(t))$ the function $\phi(t)$ vanishes on $[a,b]$ .}

Set $E(t) \equiv p_{11}(t) - p_{22}(t), \ph t \ge t_0$.

Consider the scalar Riccati equation
$$
y' + p_{12}(t) y^2 + E(t) y - p_{21}(t) = 0, \phh t\ge t_0.
$$
The solutions $y(t)$  of this equation, existing on any interval $[t_1,t_2)\ph (t_0 \le t_1 < t_2 \le + \infty)$, are connected with solutions $(\phi(t), \psi(t))$ of the system (2.3) by relations (see [3])
$$
\phi(t) = \phi(t_1) \exp\biggl\{\il{t_1}{t}[p_{12}(\tau) y(\tau) + a_{11}(\tau)]d\tau\biggr\}, \ph \phi_1(t_1) \ne 0, \ph \psi(t) = y(t) \phi(t),  \eqno (2.4)
$$
$t\in [t_1.t_2)$.

\vskip 5pt

{\bf Theorem 2.2.} {\it Let the following conditions be satisfied:

\noindent
 $p_{12}(t) \ge 0, \ph t\ge t_0$;

\noindent
$\int\limits_{t_0}^{+\infty} p_{12}(t)\exp\bigl\{-\int\limits_{t_0}^tE(\tau) d\tau\bigr\}d t = - \int\limits_{t_0}^{+\infty}p_{21}(t)\exp\bigl\{\int\limits_{t_0}^tE(\tau) d\tau\bigr\} d t = +\infty$.

\noindent
Then the system (2.3) is oscillatory.}

Proof. By analogy with the proof of Theorem 2.4 from [7] (see also [10, Corollary 3.1]).

\vskip 5pt

{\bf Theorem 2.3.} {\it Let the following conditions be satisfied:

\noindent
 $p_{12}(t) \ge 0, \ph t\in [a,b];$

\noindent
 $\int\limits_a^b \min\biggl[p_{12}(t)\exp\bigl\{-\int\limits_a^tE(\tau) d\tau\bigr\}, - p_{21}(t)\exp\bigl\{\int\limits_a^tE(\tau) d\tau\bigr\}\biggr]d t \ge \pi.$

\noindent
Then the system (2.3) is oscillatory on $[a,b]$.}

Proof. By analogy with the proof of Theorem 2.3 from [7] (see also [10, Corollary 3.2]).

In the system (1.1) substitute
$$
\phi_2 =  y_1 \phi_1, \ph \phi_3 = y_2 \phi_1, \dots, \phi_n = y_{n-1} \phi_1, \eqno (2.5)
$$
where $y_1,\dots,y_{n-1}$ are any unknown functions.
We obtain
$$
\left\{\begin{array}{l}
\phi_1' = [a_{11}(t) + a_{12}(t)y_1 + a_{13}(t) y_2 +\dots + a_{1n}(t) y_{n-1}]\phi_1,\\
\bigl(y_1' + a_{12}(t) y_1^2 + a_{13}(t) y_1 y_2 + \dots + a_{1n}(t) y_1 y_{n-1} + (a_{11}(t) - a_{22}(t)) y_1 - \\
\phantom{aaaaaaaaaaaaaaaaaaaaaaaaaaaaaaa} -a_{23}(t) y_2 - \dots -a_{2n}(t) y_{n-1} - a_{21}(t)\bigr) \phi_1 = 0,\\
\bigl(y_2' + a_{12}(t) y_2 y_1 + a_{13}(t)  y_2^2 + \dots + a_{1n}(t) y_2 y_{n-1}  - a_{32}(t) y_1 +\\
\phantom{aaaaaaaaaaaaaa}+ (a_{11}(t) - a_{33}(t)) y_2  -a_{34}(t) y_3 - \dots -a_{3n}(t) y_{n-1} - a_{31}(t)\bigr) \phi_1 = 0,\\
\centerline{\dots\dots\dots\dots\dots\dots\dots\dots\dots\dots\dots\dots\dots\dots\dots\dots\dots\dots\dots\dots\dots}\\
\bigl(y_{n-1}' + a_{12}(t) y_{n_1} y_1 + a_{13}(t) y_{n-1} y_2 + \dots + a_{1n}(t) y_{n-1}^2 + (a_{11}(t)  - \\
\phantom{aaaaa} -a_{n2}(t) y_1 a_{n3}y_2 - \dots -a_{n,n-1}(t) y_{n-2} + (a_{11}(t) - a_{nn}(t) y_{n-1} - a_{n1}(t)\bigr) \phi_1 = 0,
\end{array}
\right.
$$
It follows from here that all solutions $(y_1(t),\dots,y_{n-1}(t))$ of the system
$$
\left\{\begin{array}{l}
y_1' + a_{12}(t) y_1^2 + a_{13}(t) y_1 y_2 + \dots + a_{1n}(t) y_1 y_{n-1} + (a_{11}(t) - a_{22}(t)) y_1 - \\
\phantom{aaaaaaaaaaaaaaaaaaaaaaaaa} -a_{23}(t) y_2 - \dots -a_{2n}(t) y_{n-1} - a_{21}(t)\bigr) = 0,\\
y_2' + a_{12}(t) y_2 y_1 + a_{13}(t)  y_2^2 + \dots + a_{1n}(t) y_2 y_{n-1}  - a_{32}(t) y_1 +\\
\phantom{aaaaaaaaa}+ (a_{11}(t) - a_{33}(t)) y_2  -a_{34}(t) y_3 - \dots -a_{3n}(t) y_{n-1} - a_{31}(t) = 0,\\
\phantom{aaaaa}\dots\dots\dots\dots\dots\dots\dots\dots\dots\dots\dots\dots\dots\dots\dots\dots\dots\dots\\
y_{n-1}' + a_{12}(t) y_{n_1} y_1 + a_{13}(t) y_{n-1} y_2 + \dots + a_{1n}(t) y_{n-1}^2 + (a_{11}(t)  - \\
\phantom{a} -a_{n2}(t) y_1 a_{n3}y_2 - \dots -a_{n,n-1}(t) y_{n-2} + (a_{11}(t) - a_{nn}(t) y_{n-1} - a_{n1}(t) = 0,
\end{array}
\right. \eqno (2.6)
$$
$t \ge t_0,$ existing on any interval $[t_1,t_2)\subset [t_0,\infty)$, are connected with solutions \linebreak $(\phi_1(t),\dots,\phi_n(t))$ of the system (1.1) by relations
$$
\sist{\phi_1(t) = \phi_1(t_1)\exp\biggl\{\il{t_1}{t}[a_{11}(\tau) + a_{12}(\tau) y_1(\tau) +\dots + a_{1n}(\tau) y_{n-1}(\tau)] d \tau\biggr\},}{\phi_1(t_1) \ne 0, \ph \phi_{k+1}(t) = y_k(t) \phi_1(t), \phh t \in [t_1,t_2), \phh k =\overline{1, n-1}.} \eqno (2.7)
$$

{\bf Definition 2.3.} {\it An interval $[t_1,t_2)\subset [t_0,\infty)$ is called the maximum existence interval for a solution $(y_1(t),\dots,y_{n-1}(t))$ of the system (2.6), if $(y_1(t),\dots,y_{n-1}(t))$ exists on $[t_1,t_2)$ and cannot be continued to the right from $t_2$ as a solution of the system (2.6).
}

{\bf Lemma 2.1.} {\it Let $(y_1(t),\dots,y_{n-1}(t))$ be a solution of the system (2.6) on the finite interval $[t_1,t_2)$.  If the function $F(t) \equiv \il{t_1}{t}[a_{12}(\tau) y_1(\tau) + a_{13}(\tau) y_2(\tau) +\dots + a_{1n}(\tau) y_{n-1}(\tau)] d \tau, \linebreak t\in [t_1,t_2)$ is bounded from below on $[t_1,t_2)$, then $[t_1,t_2)$ cannot be the maximum existence interval for $(y_1(t),\dots,y_{n-1}(t))$.
}

Proof. Let $(\phi_1(t),\dots,\phi_n(t))$ be a solution of the system (1.1) with the initial values $\phi_1(t_1) = 1, \ph \psi_2(t_1) = y_1(t_1),\dots, \phi_n(t_1) = y_{n-1}(t_1)$ then by (2.5) - (2.7) $\phi_1(t) = \linebreak = \exp\biggl\{\il{t_1}{t}a_{11}(\tau) d \tau + F(t)\biggr\}, \ph t \in [t_1,t_2)$. Since $F(t)$ is bounded from below from the las equality it follows that $\phi_1(t)\ne 0, \ph t \in [t_1,t_2+\varepsilon)$ for some $\varepsilon >0$ (as far as $t_2<\infty$). Then by (2.7) and by the uniqueness theorem  $(\phi_2(t)/\phi_1(t),\dots,\phi_n(t)/\phi_1(t))$ is a solution of the system (2.6) on $[t_1,t_2+\varepsilon)$, which coincides with $(y_1(t),\dots,y_{n-1}(t))$ on $[t_1,t_2)$. It follows from here that $[t_1,t_2)$ is not the maximum existence interval for $(y_1(t),\dots,y_{n-1}(t))$. The lemma is proved.

\vskip 20pt

{\bf 3. Main results}. Let $u(t)$ and $v(t)$ be any locally integrable functions on $[t_0,\infty)$. We will say that the relation $\frac{u(t)}{v(t)}$ is well defined on $[t_0,\infty)$ if $v(t) \ne 0$ almost everywhere on $[t_0,\infty)$ and there exists a locally integrable function $\mu(t)$ on  $[t_0,\infty)$ such that $u(t) = \mu(t) v(t)$ almost everywhere on $[t_0,\infty)$. If these conditions hold then we define  $\frac{u(t)}{v(t)}$ as $\frac{u(t)}{v(t)}\stackrel {def}{=}\mu(t), \ph t \ge t_0$.

{\bf Example 3.1.} {\it If $u(t) = \sin 2t, \ph v(t) = \cos t, \ph t \ge 0$, then $\frac{u(t)}{v(t)} = 2 \sin t, \ph t \ge 0$ is well defined on $[0,\infty)$, whereas for $u(t) = \frac{1}{\sqrt{t}}, \ph u(0) = 0, \ph  v(t) = \cos t \sqrt{t}, \ph t \ge 0$ the relation $\frac{u(t)}{v(t)}$ is not well defined on $[0,\infty)$.
}

Hereafter we will always assume that the relations $\frac{a_{1,k}(t)}{a_{12}(t)}, \ph k=\overline{3,n}$ are well defined on $[t_0,\infty)$ and are absolutely continuous on $[t_0,\infty)$.

We set:
$$
A(t)\equiv a_{11}(t) - a_{22}(t) - \sum\limits_{j=3}^n\frac{a_{1j}(t)a_{j2}(t)}{a_{12}(t)},
$$

$$
B_k(t)\equiv a_{22}(t)\frac{a_{1k}(t)}{a_{12}(t)} - a_{2k}(t) - \Bigl(\frac{a_{1k}(t)}{a_{12}(t)}\Bigr)' -\sum\limits_{j=3}^n\Bigr\{\frac{a_{1j}(t)}{a_{12}(t)}\Bigl[a_{jk}(t) - a_{j2}(t)\Bigr] + a_{1j}(t)\Bigl[1 - \frac{a_{1k}(t)}{a_{12}(t)}\Bigr]\Bigr\},
$$
$$
k=\overline{3,n}, \phantom{aaaaaaaaaaaaaaaaaa}  C(t)\equiv  - a_{21}(t)- \sum\limits_{j=3}^n\frac{a_{1j}(t)a_{j_1}(t)}{a_{12}(t)}, \ph t \ge t_0.\phantom{aaaaaaaaaaaaaaa}
$$

{\bf Theorem 3.1.} {\it Assume $a_{12}(t) \ge 0, \ph t \ge t_0$. Then the following assertions are valid.

\noindent
1. If

\noindent
 A) for every $T \ge t_0$ and $\sigma_3,\dots,\sigma_{n} \in \{0,1\}$ there exist $t_2=t_2(T,\sigma_3,\dots,\sigma_{n}) > \linebreak t_1 = t_1(T,\sigma_3,\dots,\sigma_{n-1})\ge T$ such that $B_k(t)(-1)^{\sigma_k} \ge 0, \ph t\in[t_1,t_2], \ph k=\overline{3,n}.$
and

\noindent
B) $\il{t_1}{t_2}\min\biggl[ a_{12}(t)\exp\biggl\{\il{t_1}{t}-A(\tau) d\tau \biggr\}, C(t)\exp\biggl\{\il{t_1}{t} A(\tau) d\tau\biggr\}\biggr]dt \ge \pi,$

\noindent
then the system (1.1) is suboscillatory.

\noindent
2. If

\noindent
C) $B_k(t) \equiv 0, \ph t \ge t_0, \ph k=\overline{3,n}$ and

\noindent
D) $\ilp{t_0} a_{12}(t)\exp\biggl\{\il{t_1}{t}-A(\tau) d\tau \biggr\}d t =\ilp{t_0} C(t)\exp\biggl\{\il{t_1}{t} A(\tau) d\tau\biggr\}dt =\infty$

\noindent
then the system (1.1) is oscillatory.
}

Proof. In the system (1.1) substitute
$$
\phi_2 = \widetilde{\phi}_2 + \nu_3(t)\phi_3 + \nu_4(t) \phi_4+\dots+\nu_n(t)\phi_n, \eqno (3.1)
$$
where $\nu_k(t)\equiv 1 - \frac{a_{1k}(t)}{a_{12}(t)}, \ph t \ge t_0, \ph k=\overline{3,n}$.
After some simple differential and arithmetic operations we obtain
$$
\left\{\begin{array}{l} \phi_1' = a_{11}(t)\phi_1 + a_{12}(t)\widetilde{\phi}_2 + a_{12}(t)\phi_3+\dots + a_{12}(t)\phi_n,\\
\widetilde{\phi_2}' = \widetilde{a}_{21}(t)\phi_1 + \widetilde{a}_{22}(t)\widetilde{\phi}_2 + \widetilde{a}_{23}(t)\phi_3+\dots + \widetilde{a}_{2n}(t)\phi_n,\\
\phi_3' = a_{31}(t)\phi_1 + a_{32}(t)\widetilde{\phi}_2 + \widetilde{a}_{33}(t)\phi_3+\dots + \widetilde{a}_{3n}(t)\phi_n,\\
\dots\dots\dots\dots\dots\dots\dots\dots\dots\dots\dots\dots\dots\dots\dots\dots\dots\\
\phi_n' = a_{n1}(t)\phi_1 + a_{n2}(t)\widetilde{\phi}_2 + \widetilde{a}_{n3}(t)\phi_3+\dots + \widetilde{a}_{nn}(t)\phi_n, \ph t \ge t_0.
\end{array}
\right. \eqno (3.2)
$$
where
$$
\widetilde{a}_{21}(t) \equiv a_{21}(t) - \sum\limits_{j=3}^{n} \nu_j(t)a_{j1}(t), \phh \widetilde{a}_{22}(t) \equiv a_{22}(t) - \sum\limits_{j=3}^{n} \nu_j(t)a_{j2}(t),
$$
$$
\widetilde{a}_{2k}(t)\equiv a_{2k}(t) + a_{22}(t)\nu_k(t) +\Bigl(\frac{a_{1k}(t)}{a_{12}(t)}\Bigr)' - \sum\limits_{j=3}^n\nu_j(t)\Bigl[a_{jk}(t) + a_{j2}(t)\nu_k(t)\Bigr], \phh k=\overline{3,n},
$$
$$
\widetilde{a}_{jk}(t)\equiv a_{jk}(t)+ a_{j2}(t)\nu_k(t), \phh j,k = \overline{3.n}
$$

Let us prove the assertion 1.
Assume the system (1.1) is not sub oscillatory. Then by (3.1), (3.2), (2.5)-(2.7) there exists $T\ge t_0$ such that the system
$$
\left\{ \begin{array}{l}
y_1' + a_{12}(t)(y_1^2 + y_1 y_2 +\dots +y_1 y_{n-1})+ \\
\phantom{aaaaaaaaaaaaaaaaaaaa} +(a_{11}(t) - \widetilde{a}_{22}(t)) y_1 - \widehat{a}_{23}(t) y_2 -\dots - \widetilde{a}_{2n}(t) y_{n-1} - \widetilde{a}_{21}(t) = 0,\\
y_2' + a_{12}(t)(y_2 y_1 + y_2^2 +\dots +y_2 y_{n-1})- \\
\phantom{aaaaaaaaaaaaaaaaaaaa} - a_{32}(t) y_1 + (a_{11}(t) - \widetilde{a}_{33}(t)) y_2 -\dots - \widetilde{a}_{3n}(t) y_{n-1} -a_{31}(t) = 0,\\
\dots\dots\dots\dots\dots\dots\dots\dots\dots\dots\dots\dots\dots\dots\dots\dots
\dots\dots\dots\dots\dots\dots\dots\dots\dots\dots\dots\dots\\
y_{n-1}' + a_{12}(t)(y_{n-1} y_1 + y_{n-1} y_2 +\dots +y_{n-1}^2)- \\
\phantom{aaaaaaaaaaaaaaaaaaaa} - a_{n2}(t) y_1  - \widetilde{a}_{n3}(t) y_2 -\dots +(a_{11}(t) - \widetilde{a}_{nn}(t) y_{n-1} - a_{n1}(t) = 0
\end{array}
\right.
$$
has a solution $(y_1(t),\dots,y_{n-1}(t))$ on $[T,\infty)$ such that
$$
y_k(t)(-1)^{\sigma_k} \ge 0, \ph t\ge T, \ph k=\overline{2,n-1} \eqno (3.3)
$$
for sone $\sigma_k \in \{0, 1\}, \ph k=\overline{2,n-1}$. It follows from here that
$$
Y'(t) + a_{12}(t)Y^2(t) + A(t) Y(t) + C(t) + \sum\limits_{k=3}^{n}B_{k}(t) y_{k-1}(t) = 0, \ph t \ge T, \eqno (3.4)
$$
where
$$Y(t) \equiv \sum\limits_{k=1}^{n-1} y_k(t), \ph t \ge T,
$$
(since $A(t) = a_{11}(t) - \widetilde{a}_{22}(t) - a_{32}(t) -\dots - a_{n2}(t),\ph
B_3(t) = -\widetilde{a}_{23}(t) + a_{11}(t) - \widetilde{a}_{33}(t)-\dots - \widetilde{a}_{n3}(t) - A(t),\dots, B_n(t)= -\widetilde{a}_{2n}(t) - \widetilde{a}_{3n}(t)-\dots + a_{11}(t)-\widetilde{a}_{nn}(t) - A(t), \ph C(t) = - \widetilde{a}_{21}(t) - \widetilde{a}_{31}(t)-\dots - \widetilde{a}_{n1}(t)
$).
By the condition A) of the theorem there exist $t_2 > t_1 \ge T$ such that
$$
\sum\limits_{k=2}^{n-1}B_{k+1}(t)y_k(t) \ge 0, \phh t\in[t_1,t_2]. \eqno (3.5)
$$
Consider the Riccati equations
$$
u' + a_{12}(t) u^2 + A(t) u + C(t) = 0, \ph t\ge T, \eqno (3.6)
$$
$$
u' + a_{12}(t) u^2 + A(t) u + C(t) + \sum\limits_{k=2}^{n-1}B_{k+1}(t) y_k(t)  = 0, \ph t\ge T, \eqno (3.7)
$$
It follows from (3.4) that $Y(t)$ is a solution of Eq. (3.7) on $[T,\infty)$. Then applying Theorem~ 2.1 to the pair of equations (3.6) and (3.7) and taking into account (3.5) we conclude that Eq. (3.6) has a solution  $u(t)$ on $[t_1,t_2]$. Then by (2.4) the system
$$
\sist{\phi' = A(t)(t)\phi + a_{12}(t)\psi,}{\psi' = -C(t)\phi, \phh t \ge T} \eqno (3.8)
$$
is not oscillatory on $[t_1,t_2]$. On the other hand by Theorem 2.3 from the condition B) of the theorem it follows that the system (3.8) is oscillatory on $[t_1,t_2]$. The obtained contradiction completes the proof of the first assertion. Let us prove the second one. Under the restriction C) the equality (3.4) becomes
$$
Y'(t) + a_{12}(t)Y^2(t) + A(t) Y(t) + C(t)  = 0, \ph t \ge T.
$$
Therefore $Y(t)$ is a solution of Eq. (3,6) on $[T,\infty)$. Hence, by (2.4) the system (3.8) is not oscillatory. On the other hand by Theorem 2.2 from the conditions  D) it follows that the system (3.8) is oscillatory. We have obtained a contradiction, completing the proof of the theorem.

\vskip 10pt

{\bf Definition 3.1.} {\it The system (1.1) is called Lyapunov stable (asymptotic stable), if its all solutions are bounded on $[t_0,\infty)$ (vanish at $\infty$).
}

\vskip 10pt

{\bf Theorem 3.2.} {\it If $a_{jk}(t) \ge 0, \ph t \ge t_0, \ph j\ne k, \ph j,k=\overline{1,n}$, then every solution $(\phi_1(t),\dots,\phi_n(t))$  of the system (1.1) with $\phi_k(t_0) > 0, \ph k=\overline{1,n}$ satisfies the relations
$$
\phi_k(t)\ge \phi_k(t_0)\exp\biggl\{\il{t_0}{t}a_{kk}(\tau)d\tau\biggr\}, \phh t \ge t_0, \phh k=\overline{1,n}. \eqno (3.9)
$$
Therefore, the system (1.1) is  non oscillatory and in order that the system (1.1) was  Lyapumov stable (asymptotic stable)
it is necessary that
$$
\sup\limits_{t\ge t_0}\il{t_0}{t}a_{kk}(\tau)d\tau <\infty \phh \biggl(\ilp{t_0}a_{kk}(\tau)d\tau = -\infty\biggr), \ph k=\overline{1,n}. \eqno (3.10)
$$
}

Proof. The second part of the theorem follows immediately from (3.9) and (3.10). Let us prove (3.9). Let $(y_1(t),\dots,y_{n-1}(t))$ be a solution of the system (2.6) with
$$
y_k(t_0) = \frac{\phi_{k+1}(t_0)}{\phi_1(t_0)} > 0, \phh k =\overline{1,n-1}.
$$
Show that $(y_1(t),\dots,y_{n-1}(t))$ exists on $[t_0,\infty)$ and
$$
y_k(t) > 0, \ph t \ge t_0, \ph k=\overline{1.n-1}. \eqno (3.11)
$$
Let $[t_0,T)$ be the maximum existence interval for $(y_1(t),\dots,y_{n_1}(t))$. Show that
$$
y_k(t) > 0, \phh t  \in [t_0,T), \phh k=\overline{1,n-1}. \eqno (3.12)
$$
Suppose this is not so. Then since by the condition of the theorem $y_k(t_0) > 0, \ph k=\overline{1,n}$ there exist
$k_0\in \{1,2,\dots,n-1\}$ and $t_1\in (t_0,T)$ such that
$$
y_k(t) > 0, \ph t\in [t_0,t_1), \ph k=\overline{1,n-1}, \eqno (3.13)
$$
$$
y_{k_0}(t_1) = 0. \eqno (3.14)
$$
By (2.6) we have

$$
y_{k_0}'(t) + a_{1,k_0+1}(t) y_{k_0}^2(t) + D(t) y_{k_0}(t)- \sum\limits_{j=1,j\ne k_0}^{n-1} a_{k_0+1,j+1}(t) y_j(t) - a_{k_0+1,1}(t) = 0, \ph t \in [t_0,t_1],
$$
where $D(t)\equiv \sum\limits_{j=1, j\ne k_0}^{n-1} a_{1,j+1}(t) y_j(t) + a_{11}(t) - a_{k_0+1,k_0,+1}(t), \ph t \in [t_0,t_1]$.
Consider the Riccati equations
$$
 u' + a_{1,k_0+1}(t) u^2 + D(t) u- \sum\limits_{j=1,j\ne k_0}^{n-1} a_{k_0+1,j+1}(t) y_j(t) - a_{k_0+1,1}(t) = 0, \ph t \in [t_0,t_1], \eqno (3.15)
$$
$$
 u' + a_{1,k_0+1}(t) u^2 + D(t) u = 0, \ph  t \in [t_0,t_1], \eqno (3.16)
$$
Since $u_*(t)\equiv 0$ is a solution of Eq. (3.16) and  $a_{1,k_0+1}(t) \ge 0, \ph t \ge t_0$ by the  comparison Theorem 2.1 and by the uniqueness theorem the solution $u_0(t)$ of Eq. (3.16) with $u_0(t_0) = y_{k_0}(t_0)$ exists on $[t_0,t_1]$ and
$$
u(t) > 0, \phh t \in [t_1,t_2]. \eqno (3.17)
$$
By (3.13) from the conditions of the theorem it follows that
$$
\sum\limits_{j=1,j\ne k_0}^{n-1} a_{k_0+1,j+1}(t) y_j(t) + a_{k_0+1,1}(t) \ge 0, \ph t \in [t_1,t_2]. \eqno (3.18)
$$
Obviously $u_1(t)\equiv y_{k_0}(t)$ is a solution of Eq. (3.15) on $[t_0,t_1]$. Then applying Theorem 2.1 to the pair of equations (3.15), (3.16) by taking into account (3.17) and (3.18) we obtain that $y_{k_0}(t_1) > 0$, which contradicts (3.14). The obtained contradiction proves (3.12). Show that
$$
T=\infty. \eqno (3.19)
$$
Suppose $T<\infty$. Then by Lemma 2.1 from the conditions of the theorem and from (3.12) it follows that $[t_0,T)$ is not the maximum existence interval for $(y_1(t),\dots,y_{n_1}(t))$, which contradicts our supposition. The obtained contradiction proves (3.19). Then by (2.7) we have
$$
\phi_1(t) = \phi_1(t_0)\exp\biggl\{\il{t_0}{t}[a_{11}(\tau) + a_{12}(\tau) y_1(\tau)+\dots+a_{1n}(\tau)y_{n-1}(\tau)]d \tau\biggr\}, \phh t \ge t_0.
$$
By the conditions of the theorem it follows from here, from (3.12) and (3.19) that.
$$
\phi_1(t) \ge \phi_1(t_0)\exp\biggl\{\il{t_0}{t}a_{11}(\tau)d \tau\biggr\}, \phh t \ge t_0.
$$
By analogy can be proved (after interchanging  $\phi_1(t)$ and $\phi_k(t)$) that
$$
\phi_k(t) \ge \phi_k(t_0)\exp\biggl\{\il{t_0}{t}a_{kk}(\tau)d \tau\biggr\}, \phh t \ge t_0, \ph k=2,3,\dots, n.
$$
The relations (3.9) are proved and the proof of the theorem is completed.

\vskip 20 pt

\centerline{ \bf References}

\vskip 10pt

\noindent
1. K. I. Al - Dosary, H. Kh. Abdullah and D. Husein. Short note on oscillation of matrix \linebreak \phantom{a} Hamiltonian systems. Yokohama Math. J., vol. 50, 2003.

\noindent
2. Sh. Chen, Z. Zheng, Oscillation criteria of Yan type for linear Hamiltonian systems, \linebreak \phantom{a} Comput.  Math. with Appli., 46 (2003), 855--862.

\noindent
3. G. A. Grigorian, On the Stability of Systems of Two First - Order Linear Ordinary\linebreak \phantom{a} Differential Equations, Differ. Uravn., 2015, vol. 51, no. 3, pp. 283 - 292.

\noindent
4. G. A. Grigorian,
Oscillatory and non-oscillatory criteria for linear four-dimensional \linebreak \phantom{aa} Hamiltonian systems
Mathematica Bohemica, Vol. 146, No. 3, 2021,  pp. 289-304.

\noindent
5. G. A. Grigorian,   Criteria of global solvability for Riccati scalar equations. Izv. Vyssh. \linebreak \phantom{aa} Uchebn. Zaved. Mat., 2015, Number 3, Pages 35–48.

\noindent
6. G. A. Grigorian, Oscillation criteria for linear matrix Hamiltonian systems. \linebreak \phantom{aa} Proc. Amer. Math. Sci, Vol. 148, Num. 8 ,2020, pp. 3407 - 3415.

\noindent
7. G. A. Grigorian.   Interval oscillation criteria for linear matrix Hamiltonian systems,\linebreak \phantom{a}  Rocky Mount. J. Math.,  vol. 50 (2020), No. 6, 2047–2057

\noindent
8.  G. A. Grigorian. Oscillatory criteria for the systems of two first - order Linear \linebreak \phantom{a} ordinary differential equations. Rocky Mount. J. Math., vol. 47, Num. 5,
 2017, \linebreak \phantom{a}  pp. 1497 - 1524

\noindent
9. G. A. Grigorian,  On two comparison tests for second-order linear  ordinary\linebreak \phantom{aa} differential equations (Russian) Differ. Uravn. 47 (2011), no. 9, 1225 - 1240; trans-\linebreak \phantom{aa} lation in Differ. Equ. 47 (2011), no. 9 1237 - 1252, 34C10.

\noindent
10. G. A. Grigorian, Oscillatory and Non Oscillatory criteria for the systems of two \linebreak \phantom{aa}   linear first order two by two dimensional matrix ordinary differential equations. \linebreak \phantom{aa}   Arch.  Math., Tomus 54 (2018), PP. 189 - 203.

\noindent
11. I. S. Kumary and S. Umamaheswaram, Oscillation criteria for linear matrix \linebreak \phantom{aa} Hamiltonian systems, J. Differential Equ., 165, 174--198 (2000).

\noindent
12. L. Li, F. Meng and Z. Zheng, Oscillation results related to integral averaging technique\linebreak \phantom{a} for linear Hamiltonian systems, Dynamic Systems  Appli. 18 (2009), \ph \linebreak \phantom{a} pp. 725--736.

\noindent
13. Y. G. Sun, New oscillation criteria for linear matrix Hamiltonian systems. J. Math. \linebreak \phantom{a} Anal. Appl., 279 (2003) 651--658.

\noindent
14. Q. Yang, R. Mathsen and S. Zhu, Oscillation theorems for self-adjoint matrix \linebreak \phantom{a}   Hamiltonian
 systems. J. Diff. Equ., 19 (2003), pp. 306--329.

\noindent
15. Z. Zheng, Linear transformation and oscillation criteria for Hamiltonian systems. \linebreak \phantom{a} J. Math. Anal. Appl., 332 (2007) 236--245.

\noindent
16. Z. Zheng and S. Zhu, Hartman type oscillatory criteria for linear matrix Hamiltonian  \linebreak \phantom{a} systems. Dynamic  Systems  Appli., 17 (2008), pp. 85--96.

\end{document}